\documentclass[11pt,a4paper,authoryear,times]{article}
 \RequirePackage[OT1]{fontenc}
\usepackage{amsmath}
\usepackage{amssymb}
\usepackage{theorem}
\usepackage{xspace}
\usepackage{geometry}
\usepackage{bm}
\usepackage{natbib}
\usepackage{bbm}
\usepackage{color}
\usepackage[utf8]{inputenc}    

\def\pth#1{\left(#1\right)}
\def\acc#1{\left\{#1\right\}}

\def\ebo{\textrm{\mathversion{bold}$\mathbf{\beta^0}$\mathversion{normal}}}

\def\eb{\textrm{\mathversion{bold}$\mathbf{\beta}$\mathversion{normal}}}

\def\eO{\textrm{\mathversion{bold}$\mathbf{\Omega}$\mathversion{normal}}} 
\def\eR{I\!\!R}
\def\eE{I\!\!E}

\def\e1{1\!\!1}

\def\XX{\textrm{\mathversion{bold}$\mathbf{X}$\mathversion{normal}}}
\def\UU{\textrm{\mathversion{bold}$\mathbf{u}$\mathversion{normal}}}
\def\xx{\textrm{\mathversion{bold}$\mathbf{x}$\mathversion{normal}}}

\theoremstyle{plain}

\newtheorem{theorem}{Theorem}[section]

\newtheorem{lemma}{Lemma}[section]

\newtheorem{remark}{Remark}[section]
 
\numberwithin{equation}{section}
  
\def\ee1{\textrm{\mathversion{bold}$\mathbf{\varepsilon}$\mathversion{normal}}}

\def\X{{\bf {X}}}
\def\XX{\textrm{\mathversion{bold}$\mathbf{X}$\mathversion{normal}}}

\def\J{{\bf J}}
\def\eR{{\bf R}}
\def\er{{\bf r}}

\def\S{{\bf {S}}}

\def\eg{\overset{.}{\textbf{g}}}
\def\egg{\overset{..}{\textbf{g}}}

\newcommand{\N}{\mathbb{N}}
\newcommand{\R}{\mathbb{R}}
\newcommand{\PP}{\mathbb{P}}

\def\argmin{\mathop{\mathrm{arg\,min}}} 
\def\hh{ \hspace*{0.5cm}}
\begin{document}
\title {{\bf Real time change-point detection in a nonlinear quantile model}}
\author{\textbf{Gabriela Ciuperca} \thanks{Universit\'e de Lyon, Universit\'e Lyon 1, CNRS, UMR 5208, Institut Camille Jordan, Bat.  Braconnier, 43, blvd du 11 novembre 1918, F - 69622 Villeurbanne Cedex, France.
 \newline  Email: Gabriela.Ciuperca@univ-lyon1.fr
\vspace{6pt}}
 \\  {{Universit\'e Lyon 1,  UMR 5208, Institut Camille Jordan, France}}\\
 } 
 \date{}
\maketitle

\noindent{\textbf{Abstract:}} Most studies in  real time change-point detection either focus on the linear model or use the CUSUM method under classical assumptions on model errors. This paper considers the sequential change-point detection in a nonlinear quantile model. A test statistic based on the CUSUM of the quantile process subgradient is proposed and studied. Under null hypothesis that the model does not change, the asymptotic distribution of the test statistic is determined. Under alternative hypothesis that at some unknown observation there is a change in model, the proposed test statistic converges in probability to $\infty$. These results allow to build the critical regions on open-end and on closed-end procedures. Simulation results, using Monte Carlo technique, investigate the performance of the test statistic, specially for  heavy-tailed error distributions. We also compare it with the classical CUSUM test statistic.  \\

\section{Introduction}
Models with heavy-tailed errors are commonly encountered  in applications. To address this problem, a very interesting possibility is the quantile method, which has as a particular case, the least absolute deviation (median) method. Quantile framework is also useful for regressions when the model errors don't meet the classical conditions:  zero mean and bounded variance. On the other hand, in applications, it is possible that the model changes to unknown observations. We so obtain  a change-point model.  There are two types of change-point problem: a posteriori and a priori (sequential). \\
A posteriori change-point problem arises when the data are completely known at the end of the experiment to process. \\
In the sequential change-point problem, which will be presented here, the change detection is performed in real time. In order to detect a possible change in model under classical suppositions on the errors, the most used technique   is the CUSUM method. \\
\hh For the two types of change-point models, the number of publications in the last years is very extensive. Given the contribution of the present paper, we prefer to  only mention references concerning either the change-point detection in the quantile models or the sequential detection of a change-point in a model.
\cite{Oka.Qu.11} consider the estimation of multiple structural changes in a linear regression quantile, while in \cite{Su.Xiao.08} a test statistic is proposed. Several tests are also considered by  \cite{Qu.08}, \cite{Furno.12}    for structural change in linear regression quantiles. Recently, new contributions to the change-point estimation in quantile linear models were made by \cite{Aue.Cheung.Lee.Zhong.14},   \cite{Zhang.Wang.Zhu.16}.  The CUSUM method can also used in  a posteriori change-point model. So, \cite{Zhang.Wang.Zhu.14} proposed a test based on the CUSUM of the subgradient of the quantile objective function for testing a posteriori presence of a change-point due to a covariate threshold. \cite{Wang.He.07} give a possible application of a quantile framework for detecting the differential expressions in GeneChip microarray studies by proposing a rank test in a posteriori change-point linear model. Another application of a posteriori change-point quantile regression for longitudinal data is given in \cite{Li.Dowling.Chappell.15} for a study of cognitive changes in Alzheimer's disease. \\
\hh The sequential detection of a change-point in a linear model based on CUSUM of the least squares (LS) residuals  was considered by \cite{Horvath.Huskova.Kokoszka.Steinebach.04}. Their results are improved later by  bootstrapping in  \cite{Huskova.Kirch.12}.  We find  the same method in \cite{Xia.Guo.Zhao.09} for a generalized linear model. The CUSUM with adaptive LASSO residuals is used by   \cite{Ciuperca.15} for  real time change-point detection in a linear model  with a large number of explanatory variables. \cite{Zou-Wang-Tang.15} propose and study a test statistic in a linear model based on subgradient of the quantile process. \\
\hh While most of the contributions to change-point model have been focused on the linear models, \cite{Ciuperca.13} considers the sequential detection of a change-point in a nonlinear model based on CUSUM of the least squares residuals. For a posteriori change-point nonlinear model,  \cite{Boldea.Hall.13} consider   the least square method, \cite{Ciuperca:Salloum:15} consider the empirical likelihood test and \cite{Ciuperca.16} the quantile method. \\
In this paper, the sequential change-point detection in a nonlinear model is studied, when the errors don't satisfy the classical conditions. \\
\hh The rest of this article is organized as follows. In Section 2, we introduce the model assumptions and the nonlinear model under the null and alternative hypotheses. We propose and study the asymptotic behavior of the test statistic. The test statistic is the CUSUM of the subgradient of the quantile objective function. For detecting the change location in model, we propose a stopping time from which the null hypothesis is rejected. In Section 3, simulation results illustrate the performance of the proposed test statistic. Two lemmas and the  proofs of the main results are given in Section 4.

 \section{Models and main results}
 In this section,   we propose and study a test statistic for detecting in real time a change in a nonlinear quantile model. The asymptotic distribution of the test statistic under null hypothesis will allow to build the asymptotic critical region.  Notations and assumptions are also given. \\
 \hh We begin by some general notations. Throughout the paper, $C$ denotes a positive generic constant not dependent on $m$, which may take different values in different formulas or even in different parts of the same formula. All vectors and matrices denoted by bold symbols and all vectors written column-wise. 
For a  $r$-vector $\textbf{v}=(v_1, \cdots , v_r)$, $\|\textbf{v} \|$ is its  euclidean norm, $\|\textbf{v} \|_1$ is its   $L_1$-norm and  $\|\textbf{v} \|_{\infty}=\max(|v_1|, \cdots , |v_r|)$. For a matrix $\textbf{M}$, we denote by $\| \textbf{M}\|_1$ the subordinate norm to the vector norm $\|. \|_1$. \\ 

Let us consider the following nonlinear parametric   model with independent observations
\[
Y_i=g(\X_i;\eb_{i})+\varepsilon_i, \qquad i=1,\cdots, m, \cdots, m+T_m.
\]
For observation $i$, $Y_i$ denotes the response variable, $\X_i$ is a $q$-vector of explicative variables and $\varepsilon_i$ is (unobserved) model error. For simplicity, we suppose that the regressors $\X_i$ are non random, although the results will typically hold for random $\X_i$'s independent of the $\varepsilon_i$'s and if $\X_i$ independent of $\X_j$ for $i \neq j$.  If $\X_i$ is random, then, the conditional distribution with respect to $\X_i$ must be considered.\\
The regression function $g:\Upsilon  \times {\cal B} \rightarrow \R$ is known up to the parameters  $\eb_i$, with $\eb_i \in {\cal B} \subseteq \R^p$, $\Upsilon \subseteq \R^q$. We suppose that the set  ${\cal B} $ is compact (see \cite{Koenker.05}).\\
\hh We suppose that on the first $m$ observations, no change in the model parameter has occur:
\[
\eb_{i}=\ebo, \qquad \textrm{for } i=1, \cdots, m, 
\]
with $\ebo$ the true value of the parameter on the observations $i=1, \cdots, m$. The observations $1, \cdots ,m $ are called historical observations and $(Y_i,\X_i)_{1 \leqslant i\leqslant m}$ are historical data. \\
We test   null hypothesis that, at each observation $i \geq m+1$ we have the same model as the $m$ first observations: 
\begin{equation}
\label{eq3}
H_0: \eb_{i}=\ebo, \qquad m+1 \leq i \leq m+T_m,
\end{equation}
against   alternative hypothesis that, at some (unknown) observation  $k^0_m$ there is a change, called also change-point, in model: 
\begin{equation}
 \label{eq4}
H_1: \exists k^0_m \geq 1 \quad  \textrm{such that  }
 \left\{
 \begin{array}{lll}
 \eb_{i}=\ebo, & \textrm{for } & m+1 \leq i \leq m+k^0_m \\
 \eb_{i}=\eb^1 \neq \ebo, & \textrm{for } & m+k^0_m+1  <  i \leq m+T_m.
 \end{array}
 \right.
 \end{equation}
 The values of $\eb^0, \eb^1$ are unknown and the  parameter $\eb^1$ can depend on $m$.  \\
 For the sample size, two cases are possible, which will give different results,  for the test statistics under the null hypothesis:
 \begin{itemize}
\item \textit{the open-end procedure}: $T_m= \infty$;
\item \textit{the closed-end procedure}: $T_m < \infty$, $\lim_{m \rightarrow \infty} T_m=\infty$, with $\lim_{m \rightarrow \infty} \frac{T_m}{m}=T>0$, with the possibility $T=\infty$. 
\end{itemize}
 \hh We introduce now the quantile model. For a fixed quantile index $\tau \in (0,1)$, the $\tau$th conditional quantile regression of $Y$, given $\xx$, is $g(\xx,\eb)+F^{-1}(\tau)$, with $F^{-1}(\tau)$  the $\tau$th quantile ($F^{-1}$ is the inverse of the distribution function  $F$) of error $\varepsilon$. We suppose that $F(0)=\tau$. The density of $ \varepsilon$ is denoted by $f$. \\
For a fixed quantile index  $\tau \in (0,1)$, consider the check function $\rho:\R \rightarrow \R$ given by 
\[\rho_\tau(u)  = u \psi_\tau(u),
\]
with
\[
\psi_\tau(u)=\tau - \e1_{u \leq 0}.
\]
\hh The quantile  estimator of  parameter $\eb$, calculated on the historical  observations $1, \cdots, m$, is defined as the minimizer of the corresponding quantile process:
\[
\widehat \eb_m \equiv \argmin_{\eb \in {\cal B}} \sum^m_{i=1} \rho_\tau ( Y_i - g(\X_i,\eb)).
\]

We now state the assumptions on the errors, the design and on the regression function. \\
For the errors $\varepsilon_i$ we suppose that:\\
\textbf{(A1)} The errors $(\varepsilon_i)_{1 \leq i \leq n}$ are supposed independent identically distributed (i.i.d.) random variables.  We denote by  $f$ the density and by  $F$ the  distribution function of  $\varepsilon$. \\
\textbf{(A2)} $f$ is continuous, uniformly bounded away from zero and infinity and has a bounded first derivative in the neighbourhood of $\pm \bigg(g\big(\xx;\ebo+m^{-1/2}\|\UU\| \big)-g(\xx;\ebo)\bigg)$, for all $\xx \in \Upsilon$ and any  bounded $p$-vector $\UU$. \\ 

For a bounded $p$-vector $\UU \in \R^p $, let  ${\cal V }_m(\UU) \equiv  \{ \eb \in {\cal B } ; \| \eb-\ebo \| \leq m^{-1/2} \| \UU \| \}$ be the neighbourhood of  $\ebo$.\\
\hh The regression function $g(\xx,\eb)$ is supposed continuous on $\Upsilon$, twice differentiable in $\eb$, continuous twice differentiable for $\eb \in {\cal V }_m(\UU)$, for any $\UU$ bounded. In the following, for $\xx \in \Upsilon$ and $\eb \in \Gamma$ we use notation $\eg(\xx,\eb) \equiv \partial g(\xx, \eb)/ \partial \eb$ and $\egg(\xx,\eb) \equiv \partial^2 g(\xx, \eb)/ \partial \eb^2$. Moreover, for the function $g$, following assumptions are considered:\\
\textbf{(A3)} For all $\xx \in \Upsilon$, function $\eg(\xx,\eb)$ is bounded for any $\eb$ in every   neighbourhood  ${\cal V }_m(\UU)$, for bounded  $\UU$.\\

The design $(\XX_i)_{1 \leq i \leq T_m}$ and the function $g$ satisfy the following assumptions:\\
\textbf{(A4)} There exist two positive constants $C_1,C_2 >0$ and natural $m_0$, such that, for all $\eb_1, \eb_2 \in \Gamma$ and $m \geq m_0$:
$
C_1 \| \eb_1 - \eb_2 \|_2 \leq \left( m^{-1} \sum^m_{i=1} [g(\XX_i,\eb_1) -g(\XX_i,\eb_2)]^2 \right)^{1/2} \leq C_2 \| \eb_1 - \eb_2 \|_2.
$
Moreover, we have that
 $ m^{-1}   \sum^m_{i=1} \eg(\XX_i,\ebo) \eg^t(\XX_i,\ebo)$ converges, as $m \rightarrow \infty$,  to a positive definite matrix.\\
\textbf{(A5)} For any $j =1, \cdots, p$, all $p$-vector $\UU$ bounded, we have, \\  $\lim_{m \rightarrow \infty} m^{-1} \sum^m_{i=1} \sup_{\eb \in {\cal V}_m(\UU)} \|\egg_j(\X_i;\eb)\|^2 < \infty$, with $\egg_j(\X_i;\eb)\equiv \left(\frac{\partial ^2 g(\XX_i;\eb)}{\partial \beta_j \partial \beta_l} \right)_{1 \leq l \leq p}$ and $\eb=(\beta_1, \cdots , \beta_p)$.\\
\textbf{(A6)} For any  $\UU$ bounded, we have,  $\lim_{m \rightarrow \infty} m^{-1} \sum^m_{i=1} \sup_{\eb \in {\cal V}_m(\UU)} \| \eg(\X_i;\eb)\|^3 <\infty$. \\

Assumptions (A1), (A4) and  the fact that the density $f$ is bounded  are   considered also in  \cite{Koenker.05}. They are necessary to have the consistency and asymptotic normality of $\widehat \eb_m $.\\
In the linear case, assumption (A6) becomes $\sum^m_{i=1}\| \X_i\|^3 =O(m)$, assumption also considered by    \cite{Koenker.Portnoy.87} and by \cite{Zou-Wang-Tang.15} for linear quantile regression. Assumption (A2) is considered by \cite{Zou-Wang-Tang.15} for real time detection of a change in a linear quantile model, when $f'$ is bounded in the neighbourhood of $\X_i^t \ebo$. Obviously, assumption (A5)  is needed only in nonlinear models, so, it is   for example in multiple structural change-point nonlinear quantile model, considered by \cite{Ciuperca.16}, but in a stronger assumption: for all $   \eb \in {\cal B}$, $  \xx \in \Upsilon $, it is imposed that $\| \egg(\xx; \eb)\|_1$ is bounded. \\

For a vector $\eb \in {\cal B}$, let be the following $p$-square matrix: 
\[
{\J}_m(\eb) \equiv \tau (1- \tau)\frac{1}{m} \sum^m_{i=1} \eg(\X_i; \eb ) \eg^t(\X_i;\eb )
\]
and for the quantile estimator $\widehat \eb_m$, let be the following random $p$-vector:
\begin{equation}
\label{eS}
{\S} (m,k) \equiv    {\J}^{-1/2}_m(\widehat \eb_m) \sum^{m+k}_{i=m+1}\eg(\X_i;\widehat \eb_m ) \psi_\tau(Y_i - g(\X_i;\widehat \eb_m  )).
\end{equation}
The  statistic ${\S} (m,k)$ is the  similar to the cumulative sum (CUSUM) of quantile residuals  $\widehat{\varepsilon}_i \equiv Y_i - g(\X_i;\widehat \eb_m ) $, defined by $\sum^{m+k}_{i=m+1} \widehat{\varepsilon}_i $. More precisely, statistic (\ref{eS}) is the CUSUM of the subgradient of the quantile objective function.  The common principle of the CUSUM method and of the test statistic considered in present paper,  is to assume that there is the same model as for historical data,  in each observation $m+k$, for $k \geq 1$.  \\
In order to construct a test statistic for $H_0$ against $H_1$, we consider, for a given  constant $\gamma \in [0,1/2)$,   the normalisation function $z(m,k,\gamma) \equiv m^{1/2}\pth{1+\frac{k}{m}} \pth{\frac{k}{k+m}}^\gamma$ and the following random variable:
\begin{equation}
\label{eq2}
\Gamma(m,k,\gamma) \equiv \frac{\|\S (m,k)  \|_{\infty}}{z(m,k,\gamma)}.
\end{equation}
The function $z(m,k,\gamma)$, proposed by  \cite{Horvath.Huskova.Kokoszka.Steinebach.04}, is used as a boundary function. \\ 
For testing $H_0$, given by (\ref{eq3}), against $H_1$, given by (\ref{eq4}), we will consider the following test statistic:
\begin{equation}
\label{Zmg}
 Z_m(\gamma) \equiv\sup_{1 \leq k \leq T_m} \Gamma(m,k,\gamma).
\end{equation}
\begin{remark}
For linear function $g(\xx;\eb)=\xx^t \eb$, statistic ${\S} (m,k)$ becomes:\\ 
$
\big( \tau (1-\tau)m^{-1} \sum^m_{i=1} \XX_i \XX_i^t \big)^{-1/2}  \sum^{m+k}_{i=m+1} \XX_i \psi_\tau(Y_i - g(\X_i;\widehat \eb_m  ))
$
and then $ Z_m(\gamma)$ is the test statistic studied by  \cite{Zou-Wang-Tang.15} for sequential change-point  detection   in a linear quantile model.
\end{remark}

In order to study the asymptotic behaviour of statistic $ Z_m(\gamma)$, that depends on  quantile estimator $\widehat \eb_m$, one needs  a similar expression to the Bahadur representation for the linear quantile model. For this, by assumption (A4), the following matrix can be considered:
\begin{equation}
\label{om}
\eO \equiv f(0) \lim_{m \rightarrow \infty} \frac{1}{m} \sum^m_{i=1} \eg(\X_i; \ebo ) \eg^t(\X_i; \ebo ).
\end{equation}

On the other hand, assumption (A3) implies that $\max_{1 \leq i \leq m} m^{-1/2} \| \eg(\X_i;\ebo)\|\rightarrow 0$ as $m \rightarrow \infty$. This implies  the uniqueness of $\widehat{\eb}_m$ and the  following representation
 (see   \cite{Koenker.05}, page 100) for $\widehat \eb_m $:
\begin{equation}
\label{eq1}
\widehat \eb_m = \ebo+\frac{1}{m} \eO^{-1}  \sum^m_{l=1} \eg^t(\X_l; \ebo )\psi_\tau(Y_l- g(\X_l;\ebo )) +o_{\PP}(m^{-1/2}).
 \end{equation}

This representation implies that, in order to show results where $\widehat \eb_m$ appears, we can  use instead of  estimator $\widehat \eb_m $, some parameter $\ebo+m^{-1/2} \UU$, with $\UU$ a bounded vector of size $p$. 

\begin{remark}
For real time detection in a nonlinear model by CUSUM method with LS residuals, considered by \cite{Ciuperca.13}, we must calculate the LS estimators of the model parameters. Numerically, it can cause problems, since the algorithm should start with a point in a neighbourhood  of the unknown true value. If the algorithm starts with a bad starting point, than it can not converge. For the test statistic proposed in the present paper, the quantile estimator and then the test statistic $Z_m(\gamma)$ can always be calculated.
\end{remark}

 The following theorem establishes the limiting distribution of the test statistic under   null hypothesis. This distribution depends on the value of $\gamma \in [0,1/2)$. Will be studied by numerical simulations, in Section \ref{simus}, the most appropriate values of $\gamma$. The proof of Theorem \ref{theorem 1} is given in Subsection \ref{proofs_Th}.
 
 \begin{theorem}
 \label{theorem 1}
 If hypothesis  $H_0$ holds, under assumptions (A1)-(A6), for all constant $\gamma \in [0, 1/2)$,\\
 (i)  if $T_m=\infty$ or ($T_m <\infty$ and $\lim_{m \rightarrow \infty} T_m/m=\infty$), then,  
 \[
 Z_m(\gamma)   \overset{{\cal{L}}} {\underset{m \rightarrow \infty}{\longrightarrow}}  \sup_{0< t <1} \frac{\| \textbf{W}(t) \|_{\infty}}{t^\gamma}.
 \]
(ii) if $T_m <\infty$ and $\lim_{m \rightarrow \infty} T_m/m=T< \infty$, then $ Z_m(\gamma)$ converges, as $m \rightarrow \infty$, to $ \sup_{0< t <T/(1+T)} \frac{\| \textbf{W}(t) \|_{\infty}}{t^\gamma}$.\\
 For the two statements, $\{ \textbf{W}(t), t \in [0,\infty) \}$ is a Wiener process of  dimension $p$.
 \end{theorem}
 
 \begin{remark} Beside the linear quantile model considered by \cite{Zou-Wang-Tang.15}, in our case, in order to obtain relations (\ref{eq18}) and (\ref{eq19}) in the proof of Theorem \ref{theorem 1}, one must use the   K-M-T approximation and not  relation (3) of \cite{Leisch.Hornik.13}  which gives a version of FCLT and not uniform approximations for the two sums of relations (\ref{eq15}) and (\ref{eq16}).
 \end{remark}
  
  In order to have a test statistic, thus, to build a critical region, it is necessary to study the behaviour  of $Z_m(\gamma)$ under   alternative hypothesis $H_1$. 
 For this, we suppose that, under hypothesis $H_1$, the (unknown) change-point $k^0_m$ is not very far from the last observation of historical observations. Obviously, this supposition poses no problem for practical applications, since if hypothesis $H_0$ was not rejected until an observation $k_m$ of order $m^s$, we reconsider as historical periode, all observations of 1 to $k_m$.  \\
    Consider then  the following assumption on $k^0_m$ of (\ref{eq4}).\\
\textbf{(A7)}      $k^0_m=O(m^s)$, with the constant $s$:
$\left\{  \begin{array}{lll}
s\geq 0, & & \textrm{for open-end procedure, } \\ 
0 \leq s \leq 1, & & \textrm{for closed-end procedure}.
\end{array}
\right.$\\

Obviously, an  identifiability assumption   is also necessary, that is,  if the parameters before and after the change are different, then the model is also different.  On the other hand, the jump in the parameters must be strictly greater than the  convergence rate of the quantile estimator $\widehat{\eb}_m$. Then, the following assumption is considered on the regression function $g$ and on the parameter $\eb^1$ after the change-point $k^0_m$. \\
  \textbf{(A8)}  \textit{(a)} If $\| \eb^1-\ebo\| >c_1$, with $c_1$ a positive constant not depending on $m$,  there exists a positive constant $c_2$ such that 
\[
\frac{1}{a_m} \left\| \sum ^{m+k^0_m+a_m}_{i=m+k^0_m+1}  \eg(\X_i; \ebo ) \big[ g(\X_i;\eb^1) - g(\X_i;\ebo)\big]\right\|_{\infty} >c_2>0,
\]  
for any sequence $(a_m)_{m \in \N}$ converging to infinity as $m \rightarrow \infty$.\\
\hh \hh  \textit{(b)} If $\eb^1-\ebo$ converges to $\textbf{0}$ with the rate $b_m$, furthermore $m^{1/2}b_m  \rightarrow \infty$, then we suppose that $\eg(\xx,\eb)$ is bounded for any $\xx \in \Upsilon$ and $\eb$ such that $\| \eb -\ebo \| \leq b_m$.\\
  
  Assumption (A7) is very natural and typically used in sequential change-point detection works (see \cite{Ciuperca.13}, \cite{Huskova.Kirch.12}). 
 Concerning assumption (A8) for linear models,  assumption \textit{{(b)}} becomes assumption (A3), while assumption \textit{(a)}  is true for all $\eb^1$ such that $\| \eb^1-\ebo\| >c_1>0$ by assumption (A4).  For nonlinear random models under classical assumptions on errors, in  papers of \cite{Ciuperca.13} and of \cite{Boldea.Hall.13}, where sequential and a posteriori  change-point detection  is considered, respectively, the imposed identifiability assumption is $\eE[g(\XX;\ebo)] \neq \eE[g(\XX, \eb^1)]$, for  $\ebo \neq \eb^1$.\\
  
  By the following theorem, we prove that under alternative hypothesis $H_1$, the test statistic    $Z_m(\gamma)$ converges in probability to infinity, as $m \rightarrow \infty$, for all $\gamma \in [0, 1/2)$. Its proof is given in Subsection \ref{proofs_Th}.
    
\begin{theorem}
   \label{theorem 2}
  Suppose that alternative hypothesis $H_1$ holds,   assumptions (A1)-(A6) are satisfied.   Then,  we have that:
    \[
 Z_m(\gamma)=\sup_{1 \leq k \leq T_m} \Gamma(m,k,\gamma) \overset{{\PP}} {\underset{m \rightarrow \infty}{\longrightarrow}}  \infty, \qquad \textrm{for  all } \gamma \in [0,1/2).
 \]
\end{theorem}
 
 Consequence of Theorem \ref{theorem 1} and Theorem \ref{theorem 2}, we deduce that the  statistic $Z_m(\gamma)$ can be considered for testing $H_0$ against $H_1$. So, the asymptotic critical region is $\{Z_m(\gamma) > c_\alpha(\gamma) \}$, where $c_\alpha(\gamma)$ is the $(1- \alpha)$ quantile of the distribution of $\sup_{0< t <1} \frac{\| \textbf{W}(t) \|_{\infty}}{t^\gamma}$, if $T_m=\infty$ or ($T_m <\infty$ and $\lim_{m \rightarrow \infty} T_m/m=\infty$), and of  $ \sup_{0< t <T/(1+T)} \frac{\| \textbf{W}(t) \|_{\infty}}{t^\gamma}$, if $T_m <\infty$ and $\lim_{m \rightarrow \infty} T_m/m=T< \infty$. On the other hand, for a fixed   $\alpha \in (0,1)$, the proposed test statistic has the asymptotic type I error probability (size) $\alpha$ and the asymptotic power 1.
 \begin{remark}
 The sequential detector statistic, built as the cumulative sum of the LS residuals, considered by  \cite{Ciuperca.13}, always in a nonlinear model, has an asymptotic distribution depending on regression function $g$ and on  true value $\ebo$  of the parameter  before the change-point. Then, for this test statistic, one must calculate the value of critical value for each function  $g$ and value $\ebo$. In the present work case, the asymptotic distribution of the proposed test statistic, and then the asymptotic critical value, don't depend on $g$ or $\ebo$. In contrast, them must be calculated for each value of $p$. 
 \end{remark}
 
 The stopping time for the proposed test statistic is the first observation $k$ for which $\Gamma(m,k,\gamma)$ is greater than $c_\alpha(\gamma)$: 
 \begin{equation}
 \label{hatk}
 \widehat k_m \equiv \left\{
 \begin{array}{l}
 \inf \acc{1 \leq k \leq T_m, \;\; \Gamma(m,k,\gamma)  \geq c_\alpha(\gamma)}, \\
 \infty, \qquad \textrm{if } \Gamma(m,k,\gamma)  <c_\alpha(\gamma), \;\; \textrm{for every } 1 \leq k \leq T_m.
 \end{array}
 \right.
 \end{equation}
 Then, hypothesis $H_0$ is rejected in $\widehat k_m $.
 
 \section{Simulations}
 \label{simus}
 In this section, Monte Carlo simulation studies are carried out to assess the performance of the proposed test statistic and to compare it with another real time detection method for nonlinear models. All simulations were performed using the R language. The program codes can be requested from the author.  The Brownian motions are generated using the \textit{BM} function in R package \textit{'sde"}. The following R packages  were also used: package \textit{"quantreg"} for functions \textit{nlrq} (nonlinear quantile model), \textit{rq} (linear quantile model), package \textit{"MASS"} for function \textit{ginv} which calculates the generalized inverse of a matrix, package  \textit{"expm"} for function \textit{sqrtm} which computes the matrix square root of a square matrix.  \\
\hh Firstly, the critical values for a model with two parameters are given, in open-end and in closed-end procedure. In Table \ref{fractiles_o} are given the empirical quantile (critical values) for the  open-end procedure, for six values of $\gamma$: \textit{0, 0.15, 0.25, 0.35, 0.45, 0.49} and five values for the nominal  size $\alpha \in \{0.01, 0.025, 0.05, 0.10, 0.25 \}$. For the same values of $\gamma$ and $\alpha$, we give the critical values in Table \ref{fractiles_c}, for the case $T=\lim_{m \rightarrow \infty}T_m/m=5/2$.\\
  Two models are considered: a linear model in subsection \ref{simu_lin} and the growth model in Subsection \ref{simu_nlin}.   
  Everywhere we take $m=200$, $T_m=500$  and $\ebo=(1,1)$. The  design is of dimension 1 ($q=1$) and it is obtained by considering  $X_i \sim {\cal N}(0,1)$, for any ${i=1,\cdots m,  \cdots, m+T_m}$.\\
 \hh   In  view of the asymptotic distribution of the test statistic  under hypothesis $H_0$, we deduce that for a fixed nominal  size $\alpha$, we have the following  relation between the empirical powers, sizes for the  open-end and  closed-end procedures: $\widehat{\pi}_{open} \leq \widehat{\pi}_{closed}$.
  \begin{table}
 \caption{\footnotesize The empirical $(1-\alpha)$ quantiles (critical values) $c_\alpha(\gamma)$ of the random variable $\sup_{0 < t < 1} \| \textbf{W}(t)\|_{\infty}/t^\gamma$, for $p=2$,  calculated on  $50000$ Monte-Carlo replications. Open-end procedure.}
\begin{center}
{\scriptsize
\begin{tabular}{|c|ccccc|}\hline  
$\gamma \downarrow ; \alpha \rightarrow$ & 0.01 & 0.025 & 0.05 & 0.10 & 0.25 \\ \hline 
0 & 2.9963 & 2.7121  & 2.4806 & 2.2257 & 1.8277 \\
0.15& 3.0538 & 2.7921  & 2.7684 & 2.5406 & 1.9059 \\
0.25& 3.1121 & 2.8333  & 2.6103 & 2.3678 & 1.9865\\
0.35& 3.2161 & 2.9422  & 2.7273 & 2.4868 & 2.1214 \\
0.45& 3.4625 & 3.2067  & 2.9943 & 2.7675 & 2.4226 \\
0.49& 3.7316 & 3.4731  & 3.2634 & 3.0303 & 2.6801 \\ \hline 
\end{tabular} 
}
\end{center}
\label{fractiles_o} 
\end{table}

  \begin{table}
 \caption{\footnotesize The empirical $(1-\alpha)$ quantiles (critical values) $c_\alpha(\gamma)$ of the random variable  $ \sup_{0< t <T/(1+T)}  {\| \textbf{W}(t) \|_{\infty}}/{t^\gamma}$, for $p=2$,  calculated on  $50000$ Monte-Carlo replications. Closed-end procedure for $T=5/2$.}
\begin{center}
{\scriptsize
\begin{tabular}{|c|ccccc|}\hline  
$\gamma \downarrow ; \alpha \rightarrow$ & 0.01 & 0.025 & 0.05 & 0.10 & 0.25 \\ \hline 
0 & 2.5478 & 2.3058  & 2.1026 & 1.8830 & 1.5430 \\
0.15& 2.7267 & 2.4791  & 2.2665 & 2.0377 & 1.6925 \\
0.25& 2.8697 & 2.6244  & 2.4098 & 2.1743 & 1.8263\\
0.35& 3.0702 & 2.8154  & 2.5979 & 2.3646 & 2.0190 \\
0.45& 3.4099 & 3.1626  & 2.9533 & 2.7271 & 2.383 \\
0.49& 3.7263 & 3.4631  & 3.2561 & 3.0311 & 2.6740 \\ \hline 
\end{tabular} 
}
\end{center}
\label{fractiles_c} 
\end{table}
\subsection{Linear model}
\label{simu_lin}
In this subsection, we present the case already considered by \cite{Zou-Wang-Tang.15},  where only standardized Gaussian errors were considered, for the   simple linear regression function $g(x; \eb)= a+b x$, with $\eb=(a,b)$. We consider in addition of   ${\cal N}(0,1)$ law for the errors, two Cauchy distributions: ${\cal C}(0,2)$ and ${\cal C}(0,1)$. The results are given in Table \ref{lin_o_q} for the  open-end procedure and in Table \ref{lin_c_q}  for the closed-end procedure. For the    historical period, we consider  $\ebo=(1,1)$.  For the model after the change $k^0_m$, we consider two  possible values for $\eb^1$: the first possible value has been considered in the paper  \cite{Zou-Wang-Tang.15}: $\eb^1=(1,2)$,  the second considered value is  $\eb^1=(2,3)$. Compared to     \cite{Zou-Wang-Tang.15}, we give empirical powers when the change occurs in $k^0_m=5$. Concerning the empirical sizes, for $\varepsilon \sim {\cal N}(0,1)$,  we find the results given in   \cite{Zou-Wang-Tang.15}. 
\begin{table}
 \caption{\footnotesize Empirical test  sizes ($\widehat{\alpha}$) and powers ($\widehat{\pi}$),  calculated on  $5000$ Monte-Carlo replications. Open-end procedure, linear  quantile model, $k^0_m=5$, $\varepsilon \sim {\cal N}(0,1)$, $\varepsilon {\sim \cal C}(0,1)$ or $\varepsilon \sim {\cal C}(0,2)$ .}
\begin{center}
{\scriptsize
\begin{tabular}{|c|c|c|c|cccccc|}\hline  
 & & &  & \multicolumn{6}{c|}{$\gamma $} \\ 
\cline{5-10} 
$\widehat{\alpha}$ or $\widehat{\pi}$  & $\eb^1$ &  $\alpha $& $\varepsilon$  & 0 & 0.15 & 0.25 & 0.35 & 0.45 & 0.49 \\ \hline \hline
  &  &  & ${\cal N}(0,1)$ & 0.007 & 0.012 & 0.017 & 0.020 & 0.024 & 0.022 \\ 
  \cline{4-10} 
   & &  0.025 &  ${\cal C}(0,2)$ & 0.006 & 0.011 & 0.013 & 0.017 & 0.021 & 0.016 \\
   \cline{4-10} 
   & &  &  ${\cal C}(0,1)$ & 0.008 & 0.015  & 0.019  & 0.022  & 0.022  & 0.17  \\
 \cline{3-10}      
 $\widehat{\alpha}$  & $\eb^0$ & & ${\cal N}(0,1)$ & 0.017 & 0.025 & 0.032 & 0.040 & 0.041 & 0.036 \\ 
  \cline{4-10} 
 &  &  0.05 &  ${\cal C}(0,2)$ & 0.016 & 0.024 & 0.030 & 0.038 & 0.039 & 0.034 \\
   \cline{4-10} 
 &  & &  ${\cal C}(0,1)$ & 0.020 &  0.026 & 0.031  & 0.038  & 0.039  &  0.031 \\
  \hline
    &  &  & ${\cal N}(0,1)$ & 1 & 1 &1& 1 & 1 & 1 \\ 
  \cline{4-10} 
   & &  0.025 &  ${\cal C}(0,2)$ & 0.552 & 0.614 & 0.641 & 0.650 & 0.587 & 0.488 \\
   \cline{4-10} 
   & &  &  ${\cal C}(0,1)$ & 0.997 &  0.998 &  0.999 & 0.999 & 0.997  & 0.995  \\
 \cline{3-10}      
 $\widehat{\pi}$  & $(1,2)$ & & ${\cal N}(0,1)$ &  1 & 1 &1& 1 & 1 & 1 \\ 
  \cline{4-10} 
 &  &  0.05 &  ${\cal C}(0,2)$ & 0.695 & 0.738 & 0.756 & 0.755 & 0.695 & 0.593 \\
   \cline{4-10} 
 &  & &  ${\cal C}(0,1)$ & 1  & 1 & 0.1  & 1  & 0.999  & 0.997  \\
  \hline
   &  &  & ${\cal N}(0,1)$ & 1 & 1 &1& 1 & 1 & 1 \\ 
  \cline{4-10} 
   & &  0.025 &  ${\cal C}(0,2)$ & 0.994 & 0.997 & 0.997 & 0.997 & 0.995 & 0.985 \\
   \cline{4-10} 
   & &  &  ${\cal C}(0,1)$ &  1 & 1 &1& 1 & 1 & 1  \\
 \cline{3-10}      
 $\widehat{\pi}$  & $(2,3)$ & & ${\cal N}(0,1)$ &  1 & 1 &1& 1 & 1 & 1 \\ 
  \cline{4-10} 
 &  &  0.05 &  ${\cal C}(0,2)$ &  0.998 & 0.999 &1& 0.999 & 0.998 & 0.995  \\
   \cline{4-10} 
 &  & &  ${\cal C}(0,1)$ &  1 & 1 &1& 1 & 1 & 1   \\
  \hline
 \end{tabular} 
}
\end{center}
\label{lin_o_q} 
\end{table}

\begin{table}
 \caption{\footnotesize Empirical test sizes ($\widehat{\alpha}$) and powers ($\widehat{\pi}$),  calculated on  $5000$ Monte-Carlo replications. Closed-end procedure, linear quantile model, $k^0_m=5$, $\varepsilon \sim {\cal N}(0,1)$, $\varepsilon {\sim \cal C}(0,1)$ or $\varepsilon \sim {\cal C}(0,2)$ .}
\begin{center}
{\scriptsize
\begin{tabular}{|c|c|c|c|cccccc|}\hline  
 & & &  & \multicolumn{6}{c|}{$\gamma $} \\ 
\cline{5-10} 
$\widehat{\alpha}$ or $\widehat{\pi}$  & $\eb^1$ &  $\alpha $& $\varepsilon$  & 0 & 0.15 & 0.25 & 0.35 & 0.45 & 0.49 \\ \hline \hline
  &  &  & ${\cal N}(0,1)$ & 0.058 & 0.030 & 0.032 & 0.029 & 0.026 & 0.022 \\ 
  \cline{4-10} 
   & &  0.025 &  ${\cal C}(0,2)$ & 0.029 & 0.029 & 0.029 & 0.027 & 0.024 & 0.017 \\
   \cline{4-10} 
   & &  &  ${\cal C}(0,1)$ & 0.030 & 0.030  & 0.031  &  0.030  & 0.024  & 0.018  \\
 \cline{3-10}      
 $\widehat{\alpha}$  & $\eb^0$ & & ${\cal N}(0,1)$ & 0.058 & 0.057 & 0.058 & 0.056 & 0.046 & 0.037 \\ 
  \cline{4-10} 
 &  &  0.05 &  ${\cal C}(0,2)$ & 0.054 & 0.054 & 0.053 & 0.054 & 0.045 & 0.035 \\
   \cline{4-10} 
 &  & &  ${\cal C}(0,1)$ &0.055  & 0.056  & 0.051  & 0.051  &  0.045 & 0.032  \\
  \hline
    &  &  & ${\cal N}(0,1)$ & 1 & 1 &1& 1 & 1 & 1 \\ 
  \cline{4-10} 
   & &  0.025 &  ${\cal C}(0,2)$ & 0.782 & 0.767 & 0.748 & 0.715 & 0.609 & 0.493 \\
   \cline{4-10} 
   & &  &  ${\cal C}(0,1)$ & 1 &  1&  1 & 1 & 1  & 0.995  \\
 \cline{3-10}      
 $\widehat{\pi}$  & $(1,2)$ & & ${\cal N}(0,1)$ &  1 & 1 &1& 1 & 1 & 1 \\ 
  \cline{4-10} 
 &  &  0.05 &  ${\cal C}(0,2)$ & 0.860 & 0.849 & 0.833 & 0.808 & 0.715 & 0.598 \\
   \cline{4-10} 
 &  & &  ${\cal C}(0,1)$ & 1  &  1& 1  & 1  & 1  & 0.997  \\
  \hline
   &  &  & ${\cal N}(0,1)$ & 1 & 1 &1& 1 & 1 & 1 \\ 
  \cline{4-10} 
   & &  0.025 &  ${\cal C}(0,2)$ & 0.999 & 0.999 & 0.999 & 0.999 & 0.996 & 0.985 \\
   \cline{4-10} 
   & &  &  ${\cal C}(0,1)$ &  1 & 1 &1& 1 & 0.998 & 0.996  \\
 \cline{3-10}      
 $\widehat{\pi}$  & $(2,3)$ & & ${\cal N}(0,1)$ &  1 & 1 &1& 1 & 1 & 1 \\ 
  \cline{4-10} 
 &  &  0.05 &  ${\cal C}(0,2)$ &  1 & 1 &1& 0.999 & 0.999 & 0.995  \\
   \cline{4-10} 
 &  & &  ${\cal C}(0,1)$ &  1 & 1 &1& 1 & 1 & 1   \\
  \hline
 \end{tabular} 
}
\end{center}
\label{lin_c_q} 
\end{table}

Comparing the results of Tables \ref{lin_o_q} and \ref{lin_c_q}, for the two Cauchy laws ${\cal C}(0,1)$ and ${\cal C}(0,2)$, we deduct that for ${\cal C}(0,2)$, if $\eb^1=(1,2)$, "only" 495 observations were not sufficient to detect the change, when the distribution of the errors $\varepsilon$ has very heavy tails. If $\eb^1=(2,3)$, thus the difference between  $\ebo$ and $\eb^1$ is bigger, then the empirical powers are very close to 1. As argument that $T_m=500$ is not large enough to detect change when $\varepsilon \sim {\cal C}(0,2)$, we see improved results when   $T_m=1500$ observations are considered after $m=200$. We got, for the \textit{open-end procedure}, the following empirical powers:
\begin{tabbing} 
\hspace{3cm} \= \hspace{1.85cm} \= \hspace{1.85cm} \= \hspace{1.85cm} \= \hspace{1.85cm} \= \hspace{1.85cm} \=  \hspace{1.85cm} \=\kill
 \> {\footnotesize$\gamma =0$} \> {\footnotesize $\gamma =0.15$} \> {\footnotesize $\gamma =0.25$ } \> {\footnotesize $\gamma =0.35$ } \> {\footnotesize $\gamma =0.45$ }\> {\footnotesize  $\gamma =0.49$ } \\
{\footnotesize  $\alpha=0.025:$ } \> 0.856 \> 0.863  \> 0.859  \> 0.844  \> 0.774  \> 0.667 \\
 {\footnotesize  $\alpha=0.05:$ } \> 0.921 \> 0.922  \> 0.918  \> 0.903  \> 0.849  \> 0.767
\end{tabbing}
and for the \textit{closed-end procedure}:
\begin{tabbing} 
\hspace{3cm} \= \hspace{1.85cm} \= \hspace{1.85cm} \= \hspace{1.85cm} \= \hspace{1.85cm} \= \hspace{1.85cm} \=  \hspace{1.85cm} \=\kill
 \> {\footnotesize$\gamma =0$} \> {\footnotesize $\gamma =0.15$} \> {\footnotesize $\gamma =0.25$ } \> {\footnotesize $\gamma =0.35$ } \> {\footnotesize $\gamma =0.45$ }\> {\footnotesize  $\gamma =0.49$ } \\
{\footnotesize  $\alpha=0.025:$ } \> 0.952 \> 0.934  \> 0.915  \> 0.882  \> 0.796  \> 0.671 \\
{\footnotesize   $\alpha=0.05:$ } \> 0.973 \> 0.963  \> 0.954  \> 0.932  \> 0.862  \> 0.776
\end{tabbing}
For $\varepsilon \sim  {\cal N}(0,1)$ or ${ \cal C}(0,1)$, the empirical powers $\widehat{\pi}$ are either 1 or very close to  1.\\
In the  open-end  procedure (see Table \ref{lin_o_q}), for all cases, the empirical sizes $\widehat{\alpha}$ are inferior to the corresponding nominal sizes $\alpha$. In the closed-end procedure (see Table \ref{lin_c_q}), for $\gamma \in \{0, 0.15, 0.25, 0.35 \}$ the empirical sizes are slightly higher than the nominal sizes for the three laws. For $\gamma \in \{0.45, 0.49 \}$ the empirical sizes $\widehat{\alpha}$ are lower to $\alpha$.
\subsection{Nonlinear model: growth model}
\label{simu_nlin}
 Let us consider now the  growth function  $g(x;\eb)=b_1-\exp(-b_2 x)$, which models many phenomena, with the parameters  $\eb=(b_1,b_2) \in   \R^2$ and $x \in \R$. We consider the following values for the true parameters: $\ebo=(1,1)$ and $\eb^1=(1,2)$. The results  are given in Table \ref{nlin_c_q} for the closed-end procedure and in Table \ref{nlin_o_q} for the open-end procedure.\\
\hh For the open-end procedure, we compared our results with the results obtained by \cite{Ciuperca.13} using CUSUM method for the least squares (LS) residuals. Recall that for the CUSUM method for the LS residuals, in addition to $\gamma$ and $\alpha$, the critical values must be calculated for each  regression function $g(\xx,\eb)$, each value of   $\ebo$ and each law of the design  $(\XX_i)_{1 \leq i \leq m}$. The proposed test statistic $Z_m(\gamma)$ given by (\ref{Zmg}), have asymptotic distribution under hypothesis $H_0$ and then the critical values,  only depend on the value of $\gamma \in (0,1/2)$ and on the nominal size $\alpha \in (0,1)$. We will study by simulations if there are more appropriate values of $\gamma$ than other. The influence of error distribution on the results is also studied. We will consider the standard normal distribution  and the Cauchy distribution,   a law which does not have finite moments, with heavy tails. \\
\hh For the open-end procedure, the empirical sizes are inferior to  $\alpha$, for any considered value of $\gamma$, $\alpha$ and for the two error  distributions   ${\cal N}(0,1)$ and ${ \cal C}(0,1)$.\\
\hh For the closed-end procedure, we have considered three distributions for the errors: ${\cal N}(0,1)$, ${ \cal C}(0,1)$ and ${\cal C}(0,2)$.   For errors  ${ \cal C}(0,2)$, we have that  $\widehat{\alpha}$ are slightly superior to $\alpha$.  If $\varepsilon \sim {\cal N}(0,1)$ or ${ \cal C}(0,1)$, then  $\widehat{\alpha} \leq \alpha$ for $\gamma \in \{ 0.45, 0.49\}$.\\
\hh On empirical powers $\widehat{\pi}$, in the all cases, for both types of procedure, they are either 1 or very close to 1. Compared to the linear case considered in Subsection \ref{simu_lin}, the changes in model are detected almost always ($\widehat{\pi} \geq 0.988$) for  the errors ${\cal C}(0,2)$. A  possible explanation of the result improvement is the model shape. More precisely, in the nonlinear case, the derivatives $\eg(x,\eb)$ depend on  $\eb$ and therefore changes appear earlier in the test statistic. In the linear case, changes occur only in $\psi_\tau$, which lengthens  change-point detection. 
\begin{table}
 \caption{\footnotesize Empirical test sizes ($\widehat{\alpha}$) and powers ($\widehat{\pi}$),  calculated on  $5000$ Monte-Carlo replications. Closed-end procedure, growth model by quantile method, $k^0_m=5$, $\varepsilon \sim {\cal N}(0,1)$, $\varepsilon {\sim \cal C}(0,1)$, $\varepsilon {\sim \cal C}(0,2)$ .}
\begin{center}
{\scriptsize
\begin{tabular}{|c|c|c|c|cccccc|}\hline  
 & & &  & \multicolumn{6}{c|}{$\gamma $} \\ 
\cline{5-10} 
$\widehat{\alpha}$ or $\widehat{\pi}$  & $\eb^1$ &  $\alpha $& $\varepsilon$  & 0 & 0.15 & 0.25 & 0.35 & 0.45 & 0.49 \\ \hline \hline
  &  &  & ${\cal N}(0,1)$ & 0.028 & 0.028 & 0.027 & 0.026 & 0.025 & 0.022 \\ 
  \cline{4-10} 
     & & 0.025 &  ${\cal C}(0,2)$ & 0.031 & 0.032  &  0.032& 0.035  &  0.038 & 0.030  \\
     \cline{4-10} 
     & &  &  ${\cal C}(0,1)$ & 0.025 & 0.026  &  0.026& 0.023  &  0.022 & 0.017  \\
 \cline{3-10}      
 $\widehat{\alpha}$  & $\eb^0$ & & ${\cal N}(0,1)$ & 0.050 & 0.051 & 0.053 & 0.051 & 0.045 & 0.035 \\ 
  \cline{4-10} 
 &  &  0.05 &  ${\cal C}(0,2)$ & 0.057 & 0.059 & 0.059 & 0.063 & 0.057 & 0.048 \\
 \cline{4-10} 
 &  &   &  ${\cal C}(0,1)$ & 0.050 & 0.049 & 0.047 & 0.047 & 0.041 & 0.036 \\
  \hline
    &  &  & ${\cal N}(0,1)$ & 1 & 1 &1& 1 & 1 & 1 \\ 
  \cline{4-10} 
   & &  0.025 &  ${\cal C}(0,2)$ & 0.999 & 0.998 & 0.998 & 0.997 & 0.994 & 0.988 \\
   \cline{4-10} 
   & &   &  ${\cal C}(0,1)$ & 0.999 & 0.999 & 0.999 & 0.999 & 0.999 & 0.999 \\
 \cline{3-10}      
 $\widehat{\pi}$  & $(1,2)$ & & ${\cal N}(0,1)$ &  1 & 1 &1& 1 & 1 & 1 \\ 
  \cline{4-10} 
 &  &  0.05 &  ${\cal C}(0,2)$ & 0.999 & 0.999 & 0.999 & 0.999 & 0.997 & 0.994 \\
 \cline{4-10} 
 &  &    &  ${\cal C}(0,1)$ & 0.999 & 0.999 & 0.999 & 0.999 & 0.999 & 0.999 \\
  \hline
 \end{tabular} 
}
\end{center}
\label{nlin_c_q} 
\end{table}
\begin{table}
 \caption{\footnotesize Empirical test sizes ($\widehat{\alpha}$) and powers ($\widehat{\pi}$),  calculated on  $5000$ Monte-Carlo replications. Open-end procedure, growth model by quantile(Q) and LS methods, $k^0_m=5$, $\varepsilon \sim {\cal N}(0,1)$, $\varepsilon {\sim \cal C}(0,1)$ .}
\begin{center}
{\scriptsize
\begin{tabular}{|c|c|c|c|c|cccccc|}\hline  
 & & & & & \multicolumn{6}{c|}{$\gamma $} \\ 
\cline{6-11} 
$\widehat{\alpha}$ or $\widehat{\pi}$  & $\eb^1$ &  $\alpha $& $\varepsilon$ & method & 0 & 0.15 & 0.25 & 0.35 & 0.45 & 0.49 \\ \hline \hline
  &  &  & ${\cal N}(0,1)$ & Q & 0.007 & 0.01 & 0.013 & 0.017 & 0.020 & 0.019 \\ 
   &  &  &   & LS & 0.002 & 0.002 & 0.002 & 0.003 & 0  & 0  \\ 
  \cline{5-11} 
     & & 0.025 &  ${\cal C}(0,1)$ & Q & 0.006 & 0.009  &  0.013 & 0.017  &  0.019 & 0.017  \\
     & &  &    & LS & 0.3 & 0.317  &  0.35 & 0.372  &  0.395 & 0.286  \\
 \cline{3-11}      
 $\widehat{\alpha}$  & $\eb^0$ & & ${\cal N}(0,1)$ & Q &  0.015 & 0.022 & 0.028 & 0.036 & 0.039 & 0.034 \\
  &  &  &   & LS & 0.005 & 0.003 & 0.002 & 0.004 & 0  & 0.001  \\  
  \cline{5-11} 
 &  &  0.05 &  ${\cal C}(0,1)$ & Q &  0.013 & 0.021 & 0.027 & 0.032 & 0.037 & 0.033 \\
    & &  &    & LS & 0.327 & 0.352  &  0.374 & 0.404  &  0.42 & 0.30  \\
  \hline
    &  &  & ${\cal N}(0,1)$ & Q & 1 & 1 &1& 1 & 1 & 1 \\ 
     & &  &    & LS & 0.672 & 0.704  &  0.685 & 0.738  &  0.719 & 0.675  \\
  \cline{5-11} 
   & &  0.025 &  ${\cal C}(0,1)$ & Q & 999 & 0.999 & 0.999 & 0.999 & 0.999 & 0.999\\
     & &  &    & LS & 0.672 & 0.704  &  0.685 & 0.738  &  0.719 & 0.675  \\
      \cline{3-11}      
 $\widehat{\pi}$  & $(1,2)$ & & ${\cal N}(0,1)$ & Q & 1 & 1 &1& 1 & 1 & 1  \\
  &  &  &   & LS & 1 & 1 &1& 1 & 1 & 1   \\  
  \cline{5-11} 
 &  &  0.05 &  ${\cal C}(0,1)$ & Q &  999 & 0.999 & 0.999 & 0.999 & 0.999 & 0.999  \\
    & &  &    & LS & 0.701 & 0.731  &  0.714 & 0.754  &  0.74 & 0.703  \\
  \hline
 \end{tabular} 
}
\end{center}
\label{nlin_o_q} 
\end{table}
\begin{table}
 \caption{\footnotesize Descriptive statistics for the stoping times,  calculated on  $5000$ Monte-Carlo replications. Open-end procedure, growth model by quantile method, $k^0_m=5$ or $k^0_m=50$, $\varepsilon \sim {\cal N}(0,1)$, $\varepsilon {\sim \cal C}(0,1)$ .}
\begin{center}
{\scriptsize
\begin{tabular}{|c|c|c|c|cccccc|}\hline  
 & & & & \multicolumn{6}{c|}{$\gamma $} \\ 
\cline{5-10} 
$k^0_m$  &  $\alpha $& $\varepsilon$ & stat & 0 & 0.15 & 0.25 & 0.35 & 0.45 & 0.49 \\ \hline \hline
 & & & median & 97 & 77 & 65 & 54 & 47 & 48 \\
 & &${\cal N}(0,1)$ & min & 43 & 26 & 19 & 11 & 1 & 1 \\
 & & & max & 230 & 214 & 190 & 184 & 187 & 200 \\
\cline{3-10} 
 & 0.025 & & median & 124 & 100 & 85 & 71 & 63 & 66 \\
 & &${\cal C}(0,1)$ & min & 51 & 31 & 25 & 1  & 1 & 1 \\
 & & & max & 437 & 421 & 413 & 413 & 437 & Inf \\
 \cline{2-10} 
 5 & & & median & 86 & 68 & 58 & 48 & 42 & 43 \\
 & &${\cal N}(0,1)$ & min & 36 & 24 & 15 & 2 & 1 & 1 \\
 & & & max &  200 & 184 & 172 & 171 & 171 & 184 \\
\cline{3-10} 
 & 0.05 & & median & 109 & 87 & 74 & 62 & 54 & 56 \\
 & &${\cal C}(0,1)$ & min & 46 & 28 & 19 & 1 & 1 & 1 \\
 & & & max &  405 & 398 & 317 & 308 & 411 & 437 \\  \hline
  & & & median & 162 & 146 & 136 & 129 & 128 & 134 \\
 & &${\cal N}(0,1)$ & min & 85 & 66 & 61 & 37 & 1 & 1 \\
 & & & max & 316 & 296 & 293 & 290 & 296 & 330 \\
\cline{3-10} 
 & 0.025 & & median & 194 & 174 & 163 & 154 & 154 & 164 \\
 & &${\cal C}(0,1)$ & min & 82 & 45 & 37 & 1  & 1 & 1 \\
 & & & max & Inf & Inf & 469 & 469 & Inf & Inf \\
 \cline{2-10} 
 50 & & & median & 149 & 134 & 127 & 120 & 119 & 125 \\
 & &${\cal N}(0,1)$ & min & 75 & 63 & 45 & 1 & 1 & 1 \\
 & & & max &  293 & 279 & 277 & 2464 & 283 & 293 \\
\cline{3-10} 
 & 0.05 & & median & 175 & 158 & 148 & 141 & 141 & 150 \\
 & &${\cal C}(0,1)$ & min & 63 & 57 & 35 & 1 & 1 & 1 \\
 & & & max &  456 & 443 & 442 & 442 & 467 & Inf \\
   \hline
 \end{tabular} 
}
\end{center}
\label{nlin_o_estim} 
\end{table}
In Table  \ref{nlin_o_estim} are given descriptive statistics  of the detected  stopping times $\widehat{k}_m$, given by (\ref{hatk}), for two values of  $k^0_m$: 5 and 50.  The value \textit{"Inf"} for  \textit{"max"}   means that there has been at least one Monte Carlo realization for which no changes were detected under  hypothesis $H_1$ (corroborating with  the empirical power  results for error law ${\cal C}(0,1)$ in Table  \ref{nlin_o_q}). 
We observe that the shortest delay time is achieved for $\gamma=0.45$. This delay time  is longer than that obtained by CUSUM method for LS  normal errors (see  \cite{Ciuperca.13}). 

\subsection{Conclusion on simulations}
 For linear model, in the open-end procedure case, the empirical sizes $\hat \alpha$ are inferior to the nominal sizes $\alpha$. In the closed-end procedure case, for $\gamma \in \{ 0, 0.15, 0.25, 0.35 \}$, the empirical sizes are slightly than the nominal size and for $\gamma \in \{0.45, 0.49   \}$ we have $\hat \alpha \leq \alpha$. For very heavy-tailed errors, if $T_m$ is not large enough and $\| \eb^1- \ebo\|$ is small, then the proposed test statistic does not detect the change.\\
\hh In contrast, for models with derivatives   dependent on   parameter $\eb$, the empirical powers are either 1 or very close to one. All values of $\gamma$ between 0 and 1/2 can be considered, with a slight advantage for the values close to 0.49, since there are fewer false alarms for heavy-tailed errors. The preference for value $\gamma = 0.45$ is given by the shortest delay for detecting the change. \\
\hh Not least, the proposed test statistic gives very good results in comparison to the CUSUM method for the LS residuals.
\section{Lemmas and proofs of Theorems}
This section contains two lemmas and the proofs of Theorems \ref{theorem 1} and \ref{theorem 2}.
\subsection{Lemmas}
In order to prove the main results, we need two lemmas. \\

By Hoeffding's inequality (see \cite{Hoeffding.63}, \cite{Wang.He.07}) we have the following result:

\begin{lemma}
\label{Lemma 4.4.1.}
Let $(Z_{in})_{1 \leqslant i \leqslant n }$ be a sequence of bounded independent random   variables. Denote by $S_n \equiv \sum^n_{i=1} Z_{in}$. Let $(\delta_n)_{n\in \N}$ be a given sequence and  suppose that 
\[
\sum^n_{i=1} Var[Z_{in}] \leq \delta_n
, \qquad \max_{1 \leq i \leq n} |Z_{in}| \leq \frac{\log 2}{2}\frac{\sqrt{\delta_n}}{\sqrt{\log n}}. \]
Then, for all $L>0$, there exists $N_0 \in \N$, such that for any $n \geq N_0$, we have:
\[
\PP \big[|S_n- \eE[S_n]| \geq (1+L) \sqrt{\delta_n} \sqrt{\log n}\big] \leq 2 \exp(- L \log n)=2 n^{-L}.
 \]
\end{lemma}

For a $p$-vector $\UU$, observations $i=m+1, \cdots , m+T_m$, let us consider the following random $p$-vectors, which   will be used in the following Lemma and in the proofs of Theorems \ref{theorem 1} and \ref{theorem 2}.
\begin{equation}
\label{Riu}
\eR_i(\UU) \equiv \eg(\X_i;\ebo+m^{-1/2} \UU)\psi_\tau(Y_i-g(\X_i; \ebo+m^{-1/2} \UU)) - \eg(\X_i;\ebo)\psi_\tau(\varepsilon_i)
\end{equation}
and, for $k=1, \cdots , T_m$, 
\begin{equation}
\label{rmk}
\er_{m,k}(\UU) \equiv \sum^{m+k}_{i=m+1} \eR_i(\UU).
\end{equation}
 \begin{lemma}
 \label{Lemma 1}
 Suppose that assumptions (A1), (A2), (A3), (A5) are satisfied. 
Let  $C_1>0$ be some constant. For all constant $L>0$ and any $k \in \N$ sufficiently large, under hypothesis $H_0$, we have:
 \[
 \PP \big[\sup_{\| \UU\| \leq C_1} \|\er_{m,k}(\UU) -\eE[\er_{m,k}(\UU)]   \| \geq  \sqrt{p}(1+L) m^{-1/4} k^{1/2}\sqrt{\log k} \big] \leq 2 k^{-L}.
 \]
 \end{lemma}
 \textbf{Proof of Lemma \ref{Lemma 1}}\\ 
 In order to prove the lemma, we show that all conditions of Lemma \ref{Lemma 4.4.1.} are satisfied for the sequence of independent random variables $\big(R_{ij}(\UU)\big)_{m+1 \leqslant i \leqslant m+k}$, for any $j=1, \cdots , p$, with $R_{ij}(\UU)$ the $j$th coordinate of $\eR_i(\UU)$.\\
 By assumption  (A3), since the function $g  $ is continuous twice differentiable in any neighbourhood ${\cal V }_m(\UU)$ of $\ebo$, with $\UU$ a  bounded $p$-vector, we have that the random vectors  $\big(\eR_i(\UU) \big)_{m+1 \leqslant i \leqslant m+k}$ are bounded with probability one:
 \[
 \max_{m+1 \leq i \leq m+k} \bigg( \max_{\| \UU\| \leq C_1} \| \eR_i(\UU) \|\bigg)= O_{\PP}(1).
 \]
 We study now $\sum^{m+k}_{i=m+1}Var[R_{ij}(\UU)]$. 
 For this, we write $\eR_i(\UU)$, defined by (\ref{Riu}), under the form:
 $
 \eR_i(\UU)=\eg(\X_i;\ebo+m^{-1/2} \UU) \big[\e1_{\varepsilon_i \leq 0} -\e1_{\varepsilon_i \leq g(\X_i; \ebo+m^{-1/2} \UU)-g(\X_i; \ebo)})\big]+\big[ \eg(\X_i;\ebo+m^{-1/2} \UU)  -\eg(\X_i;\ebo)\big][\tau - \e1_{\varepsilon_i \leq 0}] $.\\
 For all $j \in \{1, \cdots , p \}$ and $i=m+1, \cdots , m+k$, we have:
 \[
 R_{ij}(\UU)=\frac{\partial g(\X_i;\ebo+m^{-1/2} \UU)}{\partial \beta_j}\bigg(\e1_{\varepsilon_i \leq 0} -\e1_{\varepsilon_i \leq g(\X_i; \ebo+m^{-1/2} \UU)-g(\X_i; \ebo)}\bigg) \qquad \qquad
 \]
 \begin{equation}
 \label{AB}
  +\bigg( \frac{\partial g(\X_i;\ebo+m^{-1/2} \UU)}{\partial \beta_j}  -\frac{\partial g(\X_i;\ebo)}{\partial \beta_j}\bigg)(\tau - \e1_{\varepsilon_i \leq 0}) \equiv A_{ij}(\UU)+B_{ij}(\UU) .
 \end{equation}
 Thus, we have for its variance,
 \begin{equation}
 \label{eq22}
 Var[R_{ij}(\UU)] \leq \eE[(A_{ij}(\UU)+B_{ij}(\UU))^2 ] \leq 2 \big( \eE[A_{ij}^2(\UU)]+\eE[B_{ij}^2(\UU)]\big).
 \end{equation}
 By the Taylor's expansion at $\ebo$ of  $\frac{\partial g(\X_i;\ebo+m^{-1/2} \UU)}{\partial \beta_j}$, with $\UU=(u_1, \cdots, u_p)$,  we obtain, applying also the    Cauchy-Schwarz inequality: 
 \[
 \eE[B_{ij}^2(\UU)] \leq \left[ m^{-1/2}  \sum^p_{l=1} u_l \frac{\partial ^2 g(\ebo+m^{-1/2} \widetilde {\UU})}{\partial \beta_j \partial \beta_l} \right]^2 \leq m^{-1} \| \UU\|^2 \left\|\egg_j(\ebo+m^{-1/2} \widetilde {\UU}) \right\|^2,
 \]
 with the  vector $\egg_j(\ebo+m^{-1/2} \widetilde {\UU}) \equiv \left( \frac{\partial ^2 g(\ebo+m^{-1/2} \widetilde {\UU})}{\partial \beta_j \partial \beta_l}\right)_{1 \leq l \leq p}$ and the  vector $ \widetilde {\UU} =(\widetilde{u}_1, \cdots , \widetilde{u}_p)$ such that $0 \leq | \widetilde{u}_j - u_j | \leq 1$, for any  $j=1, \cdots , p$.\\
 Using assumption (A5) and $\| \UU\| \leq C_1$, we obtain, uniformly in $\UU$ and in $j$:
 \begin{equation}
 \label{eq20}
 \sum^{m+k}_{i=m+1} \eE[B_{ij}^2(\UU)] =O(m^{-1}k).
 \end{equation}
For the term $A_{ij}(\UU)$ of relation (\ref{AB}), using the identity that, for all real numbers  $a$ and $b$, we have with probability one, that, $\e1_{\varepsilon_i \leq a} -\e1_{\varepsilon_i \leq b}=\e1_{\min(a,b) \leq \varepsilon_i \leq \max(a,b)} sgn(a-b)$, we obtain for $A_{ij}^2(\UU)$:
 \[
 \eE[A_{ij}^2(\UU)] \leq  \left(\frac{\partial g((\X_i;\ebo+m^{-1/2} \UU))}{\partial \beta_j} \right)^2 \eE \left[ \left( \e1_{\varepsilon_i \leq 0} -\e1_{\varepsilon_i \leq g(\X_i; \ebo+m^{-1/2} \UU)-g(\X_i; \ebo)}\right)^2 \right]
 \]
 \[
 = \left(\frac{\partial g((\X_i;\ebo+m^{-1/2} \UU))}{\partial \beta_j} \right)^2 \eE \left[  \e1_{\min(0,g(\X_i; \ebo+m^{-1/2} \UU)-g(\X_i; \ebo)) \leq \varepsilon_i \leq \max(0,g(\X_i; \ebo+m^{-1/2} \UU)-g(\X_i; \ebo))} \right].
 \]
 We have used the notation $sgn(.)$ for the sign function. 
 We assume, without loss of generality that,  
 \[\max\big(0,g(\X_i; \ebo+m^{-1/2} \UU)-g(\X_i; \ebo)\big)=g(\X_i; \ebo+m^{-1/2} \UU)-g(\X_i; \ebo).
 \]
  Then,  $\eE \bigg[  \e1_{\min(0,g(\X_i; \ebo+m^{-1/2} \UU)-g(\X_i; \ebo)) \leq \varepsilon_i \leq \max(0,g(\X_i; \ebo+m^{-1/2} \UU)-g(\X_i; \ebo))} \bigg]=F\big( g(\X_i; \ebo+m^{-1/2} \UU)-g(\X_i; \ebo) \big) -F(0)$. 
 Thus, by two Taylor expansions, first for $F$ and after for $g$, we have:
 \[
 \eE[A_{ij}^2(\UU)]  \leq \left(\frac{\partial g(\X_i;\ebo+m^{-1/2} \UU)}{\partial \beta_j} \right)^2 \big[g(\X_i; \ebo+m^{-1/2} \UU) -g(\X_i; \ebo) \big] f(b_i)
 \]
 \[ 
 \qquad \qquad \quad = \pth{\frac{\partial g(\X_i;\ebo+m^{-1/2} \UU)}{\partial \beta_j}}^2 \left| m^{-1/2} \UU^t \eg(\X_i; \ebo+m^{-1/2} \textbf{w}_{i,m}) \right|  f(b_i)
 \]
 \[
\qquad \qquad \qquad \leq \pth{\frac{\partial g(\X_i;\ebo+m^{-1/2} \UU)}{\partial \beta_j}}^2 m^{-1/2} \|\UU\| \|\eg(\X_i; \ebo+m^{-1/2} \textbf{w}_{i,m})\|f(b_i),
 \]
 with the scalar $b_i$ between $g(\X_i; \ebo+m^{-1/2} \UU)-g(\X_i; \ebo)$ and $0$, the $p$-vector $\textbf{w}_{i,m}$ between $\ebo+m^{-1/2} \UU$ and $\ebo$. 
 For the last inequality, we have used the Cauchy-Schwarz inequality. 
 The  vector $\textbf{w}_{i,m}$ is such that there exists a constant $C_4$, which does not depend on  $i$,   $j$,   $m$, such that,  $\| \textbf{w}_{i,m} \| \leq C_4$ for any $i$ and $\UU$.\\
 Then, taking into account assumption (A2), we get
 \[
  \sum^{m+k}_{i=m+1} \eE[A_{ij}^2(\UU)] \leq C m^{-1/2} \|\UU\|  \sum^{m+k}_{i=m+1}\pth{\frac{\partial g(\X_i;\ebo+m^{-1/2} \UU)}{\partial \beta_j}}^2 \|\eg(\X_i; \ebo+m^{-1/2}  \textbf{w}_{i,m}) \|.
 \]
 Thus, using assumption (A3) and the fact that  $\UU$ is bounded, we obtain, uniformly in $\UU$:
 \begin{equation}
 \label{eq21}
 \sum^{m+k}_{i=m+1} \eE[A_{ij}^2(\UU)]  \leq O(m^{-1/2}k).
 \end{equation}
 Taking into account relations (\ref{eq22}), (\ref{eq20}), (\ref{eq21}) together with assumption (A1), we get:
 \[
 \sum^{m+k}_{i=m+1} Var[ R_{ij}(\UU)]  =O(m^{-1/2}k).
 \]
 Thus we can  apply  Lemma \ref{Lemma 4.4.1.}  for the random variable   $R_{ij}(\UU)$, for any $j =1, \cdots , p$  and for $\delta_k=k m^{-1/2}$. Therefore, for all $L>0$, $\| \UU\| \leq C_1$, for $k \in \N$ sufficiently large, we get:
 \begin{equation}
 \label{eq2kL}
  \PP \bigg[\big[R_{ij}(\UU)- \eE[R_{ij}(\UU)]\big]^2 \geq (1+L)^2  m^{-1/2} k\log k\bigg] \leq 2 k^{-L}.
 \end{equation}
 On the other hand, we have, the following inequalities, for all $  L>0$, any observation  $i=m+1, \cdots m+k$, all vector $\UU$ bounded and  enough large $k,m$, that:
 \[
 \PP \bigg[\|\eR_i(\UU) -  \eE[\eR_i(\UU)] \|^2 \geq  p(1+L)^2  m^{-1/2} k\log k\bigg] \leq \PP \bigg[ \max_{1 \leqslant j \leqslant p} |R_{ij}(\UU)- \eE[R_{ij}(\UU)]|^2 \geq  (1+L)^2  m^{-1/2} k\log k\bigg]
 \]
 \[
 \leq 2 k^{-L}.
 \]
 For the last inequality we have used relation (\ref{eq2kL}). The lemma is proved.
 \hspace*{\fill}$\blacksquare$ \\ 
\subsection{Proofs of Theorems}
\label{proofs_Th}
Here we present the proofs of Theorem \ref{theorem 1} and Theorem \ref{theorem 2}.\\

 \textbf{Proof of  Theorem \ref{theorem 1} }\\
 Statements \textit{(i)} and \textit{(ii)} are shown at the same time.\\
\textit{Case $k$  small}. Under assumptions (A1)-(A4), we have that the quantile estimator  $\widehat \eb_m$ is consistent (see \cite{Koenker.05}). By similar arguments as in the proof of Lemma \ref{Lemma 1} and by Bienaymé-Tchebychev inequality, we have that $\psi_\tau(Y_i - g(\X_i;\widehat \eb_m  ))  \overset{ \PP} {\underset{m \rightarrow \infty}{\longrightarrow}} \textbf{0}$, for any $i \geq m+1$, under hypothesis $H_1$.   Then, taking into account assumption (A3),  we have that $\eg(\X_i;\widehat \eb_m ) \psi_\tau(Y_i - g(\X_i;\widehat \eb_m  ))  \overset{ \PP} {\underset{m \rightarrow \infty}{\longrightarrow}} \textbf{0}$ and thus, taking also into account relations (\ref{eS}) and (\ref{eq2}), we have that $\Gamma(m,k,\gamma)$ converges in probability to  $0$, as $m \rightarrow \infty$. We deduct that, the  sup of $\Gamma(m,k,\gamma)$ is not achieved for small values of $k$.\\   
 \textit{Case $k$  large}.   From Lemma \ref{Lemma 1} we have for $\er_{m,k}(\UU)$, defined by (\ref{rmk}), that:
  \begin{equation}
  \label{eq5}
  \er_{m,k}(\UU)=\eE[\er_{m,k}(\UU)]+O_{\PP} (m^{-1/4} k^{1/2}\sqrt{\log k} ).
  \end{equation}
    We will calculate $\eE[\er_{m,k}(\UU)]=\sum^{m+k}_{i=m+1} \eE[\eR_i(\UU)]$. By  Taylor expansions, we have:
  \[
  \begin{array}{cl}
  \eE[\eR_i(\UU)] & = \big[\eg(\X_i;\ebo)+ m^{-1/2}  \UU \egg(\X_i;\widetilde \eb)\big] \big[F(0)-F(g(\X_i;\ebo+m^{-1/2}  \UU )-g(\X_i;\ebo)) \big] \\
  &  =- \big[\eg(\X_i;\ebo)+ m^{-1/2}  \UU \egg(\X_i;\widetilde \eb)\big] \bigg\{\big[g(\X_i;\ebo+m^{-1/2}  \UU )-g(\X_i;\ebo)) \big] f(0)  \\
   \qquad \qquad  &    \qquad \qquad +2^{-1} \big[g(\X_i;\ebo+m^{-1/2}  \UU )-g(\X_i;\ebo)) \big]^2 f'(\widetilde{b}_i) \bigg\} \\
          & = - \big[\eg(\X_i;\ebo)+ m^{-1/2}  \UU^t \egg(\X_i;\widetilde \eb)\big] \bigg\{\big[ m^{-1/2}  \UU^t \eg(\X_i;\ebo)+ 2^{-1}m^{-1} \UU^t \egg(\X_i;\widetilde \eb) \UU \big]f(0)\\
    &  \qquad \qquad +2^{-1} \big[g(\X_i;\ebo+m^{-1/2}  \UU )-g(\X_i;\ebo)) \big]^2 f'(\widetilde{b}_i)   \bigg\},
  \end{array} 
  \]
  with the scalar $\widetilde{b}_i$ between $g(\X_i;\ebo+m^{-1/2}  \UU )-g(\X_i;\ebo))$ and $0$, the vector $\widetilde \eb$   between $\ebo$ and $\ebo+m^{-1/2}  \UU$. 
     Thus, for $k$ large enough, using assumptions (A1)-(A6), we have 
  \begin{equation}
  \label{eq6}
  \eE[\er_{m,k}(\UU)]=-m^{-1/2}  f(0)\UU \sum^{m+k}_{i=m+1} \eg^t(\X_i;\ebo) \eg(\X_i;\ebo)+o(k m^{-1/2}).
  \end{equation}
    On the other hand, using relation (\ref{eq1}), we consider for $\UU$ the following value:
  \[
 \UU =m^{1/2}( \widehat \eb_m - \ebo)=m^{-1/2} \eO^{-1}  \sum^m_{l=1} \eg^t(\X_l; \ebo )\psi_\tau(Y_l- g(\X_l;\ebo)) +o_{\PP}(1).
  \]
 Then, replacing in (\ref{eq5}), taking into account of (\ref{eq6}), we have:
  \[
  \er_{m,k}(m^{1/2}( \widehat \eb_m - \ebo) )= \big[-m^{-1}\eO^{-1}  \sum^m_{l=1} \eg^t(\X_l; \ebo )\psi_\tau(Y_l- g(\X_l;\ebo)) \big]\big[f(0) \sum^{m+k}_{i=m+1} \eg^t(\X_i;\ebo) \eg(\X_i;\ebo) \big]\]
  \[\qquad \qquad \qquad \qquad +O_{\PP} (m^{-1/4} k^{1/2}\sqrt{\log k})+o_{\PP}(k m^{-1/2}).
  \]
 On the other hand, by assumption (A4), for $k$ large enough, we have: 
  \[
f(0)  \sum^{m+k}_{i=m+1} \eg(\X_i;\ebo) \eg^t(\X_i;\ebo)  =k \eO(1+O(1)),
  \]
 with the matrix $\eO$ defined by (\ref{om}).  
   Hence
  \[
  \er_{m,k}(m^{1/2}( \widehat \eb_m - \ebo) )=- k m^{-1} \sum^m_{l=1} \eg^t(\X_l; \ebo )\psi_\tau(Y_l- g(\X_l;\ebo))  +O_{\PP} (m^{-1/4} k^{1/2}\sqrt{\log k} )+o_{\PP}(k m^{-1/2}),
  \]
 which implies, taking into account (\ref{Riu}) and (\ref{rmk}), that
 \begin{equation}
  \label{eq8}
  \begin{array}{l}
\displaystyle{  \J^{-1/2}_m( \ebo)\sum^{m+k}_{i=m+1} \eg(\X_i; \widehat \eb_m)\psi_\tau(Y_i- g(\X_i;\widehat \eb_m))=  \J^{-1/2}_m(\ebo)\sum^{m+k}_{i=m+1}\eg(\X_i; \ebo)\psi_\tau(Y_i- g(\X_i;\ebo))  
 } \\
  \displaystyle{
  -k m^{-1} \J^{-1/2}_m( \ebo)\sum^m_{l=1} \eg^t(\X_l; \ebo )\psi_\tau(Y_l- g(\X_l;\ebo))+O_{\PP} (m^{-1/4} k^{1/2}\sqrt{\log k} )+o_{\PP}(k m^{-1/2}) .}
  \end{array}
  \end{equation}
  By K-M-T approximation (see \cite{Komlos.Major.Tusnady.75}, \cite{Komlos.Major.Tusnady.76}, for each $m$ there exists  two independent $p$-dimensional Wiener processes on $[0,\infty)$, $\{\textbf{W}_{1,m}(t), t \in [0,\infty)\}$ and $\{\textbf{W}_{2,m}(t), t \in [0,\infty)\}$ such that, as $m \rightarrow\infty$:
   \begin{equation}
  \label{eq15}
  \sup_{1 \leq k < \infty} k^{- 1/\nu}\left\|   \J^{-1/2}_m(\ebo)\sum^{m+k}_{i=m+1}\eg(\X_i; \ebo)\psi_\tau(Y_i- g(\X_i;\ebo)) - \textbf{W}_{1,m}(k/m) \right\|_{\infty}=O_{\PP}(1)
  \end{equation}
  and 
  \begin{equation}
  \label{eq16}
    \left\|   \J^{-1/2}_m( \ebo)\sum^m_{l=1} \eg^t(\X_l; \ebo )\psi_\tau(Y_l- g(\X_l;\ebo)) -\textbf{W}_{2,m}(1) \right\|_{\infty}=o_{\PP}(m^{1/\nu}).
  \end{equation}
   Then, by relations (\ref{eq8}), (\ref{eq15}), (\ref{eq16}) and Lemma 5.3 of Horvath et al.(2004), we have:
  \begin{equation}
  \label{eq17}
  \sup_{1 \leq k <\infty} \frac{\left\|   \J^{-1/2}_m( \ebo)\displaystyle{\sum^{m+k}_{i=m+1}} \eg(\X_i; \widehat \eb_m)\psi_\tau(Y_i- g(\X_i;\widehat \eb_m)) - \big[\textbf{W}_{1,m}(k/m) - k/m\textbf{W}_{2,m}(1)  \big] \right\|_{\infty}}{z(m,k,\gamma)}  =o_{\PP}(1).
  \end{equation}
  By a similar way as in the proof of Theorem 2.1 of \cite{Horvath.Huskova.Kokoszka.Steinebach.04}, we have that (see also the proof in the linear case, Theorem 1 of  \cite{Zou-Wang-Tang.15}):\\
 \textit{(i) if $T_m=\infty$ or ($T_m <\infty$ and $\lim_{m \rightarrow \infty} T_m/m=\infty$)}, then, 
  \begin{equation}
  \label{eq18}
   \sup_{1 \leq k <T_m} \frac{\left\|\textbf{W}_{1,m}(k/m) - k/m\textbf{W}_{2,m}(1) \right\|_{\infty}}{z(m,k,\gamma)} \overset{{\cal{L}}} {\underset{m \rightarrow \infty}{\longrightarrow}} \sup_{0 \leq  t < \infty } \frac{\|\textbf{W}_1(t)-t\textbf{W}_2(1) \|_{\infty}}{(1+t)(t/(1+t))^\gamma}\overset{{\cal{L}}} {=} \sup_{0 \leq  t  \leq 1} \frac{\| \textbf{W}(t) \|_{\infty}}{t^\gamma}, 
  \end{equation}
   \textit{(ii)  if $T_m <\infty$ and $\lim_{m \rightarrow \infty} T_m/m=T< \infty$}, then, 
  \begin{equation}
  \label{eq19}
   \sup_{1 \leq k <T_m} \frac{\left\|\textbf{W}_{1,m}(k/m) - k/m\textbf{W}_{2,m}(1) \right\|_{\infty}}{z(m,k,\gamma)} \overset{{\cal{L}}} {\underset{m \rightarrow \infty}{\longrightarrow}} \sup_{0 \leq  t \leq T } \frac{\|\textbf{W}_1(t)-t\textbf{W}_2(1) \|_{\infty}}{(1+t)(t/(1+t))^\gamma}\overset{{\cal{L}}} {=} \sup_{0 \leq  t  \leq T/(1+T)} \frac{\| \textbf{W}(t) \|_{\infty}}{t^\gamma}, 
  \end{equation}
 with $\{\textbf{W}_{1}(t), t \in [0,\infty)\}$ and $\{\textbf{W}_{2}(t), t \in [0,\infty)\}$ two independent $p$-dimensional Wiener process on $[0,\infty)$.\\
 On the other hand, by Lemma 2 of \cite{Zou-Wang-Tang.15}, we have that for all $\gamma \in [0,1/2)$ and $\nu >2$:
 \[
\lim_{m \rightarrow \infty} \frac{k^{1/\nu}+k m^{1/\nu -1}+o(k m^{-1/2} )+O(m^{-1/4} k^{1/2} \sqrt{\log k})}{z(m,k,\gamma)}= 0.
 \]
   The theorem follows by taking into account  relations (\ref{eq17}),  (\ref{eq18}),   (\ref{eq19}) and by the fact that by assumptions (A4) and (A5) we have the following decomposition:
  \[
     {\J}_m(\widehat \eb_m )  \equiv \tau(1-\tau) m^{-1}\sum^m_{l=1} \eg(\X_l; \widehat \eb_m ) \eg^t(\X_l;\widehat \eb_m ) = \J_m(\ebo) + O_{\PP}(m^{-1/2}).
  \]
    \hspace*{\fill}$\blacksquare$ \\ 
    
  \textbf{Proof of  Theorem \ref{theorem 2} }\\
 The theorem is shown if we prove that  there exists an observation $k$ depending on $m$, denoted $\widetilde k$, for which we have convergence in probability of $\Gamma(m,\widetilde k,\gamma)$, as $m \rightarrow \infty$,   to infinity. \\
{\textit{(i) Open-end procedure case.}} We suppose, without loss of generality that   $k^0_m \leq{m^s}/{2}$, with $s >1$.  We  consider $\widetilde k \equiv k^0_m+m^s$. Then, for the statistic defined by (\ref{eS}), we have:
  \begin{equation}
  \label{eq9}
{\S} (m,\widetilde k) \equiv    {\J}^{-1/2}_m(\widehat \eb_m) \bigg[\sum^{m+k^0_m}_{i=m+1}\eg(\X_i;\widehat \eb_m ) \psi_\tau(Y_i - g(\X_i;\widehat \eb_m  )) +\sum^{m+\widetilde k}_{i=m+k^0_m+1}\eg(\X_i;\widehat \eb_m ) \psi_\tau(Y_i - g(\X_i;\widehat \eb_m  )) \bigg].
\end{equation}
Similar as for Theorem \ref{theorem 1}, we have that there exists a  constant $0< C < \infty$, such that:
\[
\frac{  \| {\J}^{-1/2}_m(\widehat \eb_m)\sum^{m+k^0_m}_{i=m+1}\eg(\X_i;\widehat \eb_m ) \psi_\tau(Y_i - g(\X_i;\widehat \eb_m  ))  \|_{\infty}}{z(m,k^0_m,\gamma)} \leq C,
\]
with probability converging to  1, as $m \rightarrow \infty$.\\ Since the function $h(x) = (1+x) \big(\frac{x}{1+x} \big)^\gamma$,  is increasing in $x>0$, we have:
\begin{equation}
\label{TJ}
\frac{  \| {\J}^{-1/2}_m(\widehat \eb_m)\sum^{m+k^0_m}_{i=m+1}\eg(\X_i;\widehat \eb_m ) \psi_\tau(Y_i - g(\X_i;\widehat \eb_m  ))  \|_{\infty}}{z(m,\widetilde k,\gamma)} \leq C, 
\end{equation}
with probability converging to  1, as $m \rightarrow \infty$.\\
We will study $\sum^{m+\widetilde k}_{i=m+k^0_m+1} \eg(\X_i;\widehat \eb_m ) \psi_\tau(Y_i - g(\X_i;\widehat \eb_m  ))$ of relation (\ref{eq9}). For this, we consider,    for any $p$-vector   $\UU=O(1)$,  the following random  process:
\[
\sum^{m+\widetilde k}_{i=m+k^0_m+1} \eR_i(\UU) \equiv  \sum^{m+\widetilde k}_{i=m+k^0_m+1}  \bigg[\eg(\X_i;\ebo+m^{-1/2} \UU)\psi_\tau(Y_i-g(\X_i; \ebo+m^{-1/2} \UU)) - \eg(\X_i;\ebo)\psi_\tau(\varepsilon_i)  \bigg] .
\]
Since for $i=m+k^0_m+1, \cdots , m+\widetilde k$, hypothesis $H_1$ is true, we have that:
\begin{equation}
\label{eq10}
\sum^{m+\widetilde k}_{i=m+k^0_m+1} \eR_i(\UU)  =   \sum^{m+\widetilde k}_{i=m+k^0_m+1}  \big[\eg(\X_i;\ebo+m^{-1/2} \UU)\psi_\tau(Y_i-g(\X_i; \ebo+m^{-1/2} \UU)) - \eg(\X_i;\ebo)\psi_\tau(Y_i-g(\X_i; \eb^1))  \big],
\end{equation}
with $\eR_i(\UU)$ defined by (\ref{Riu}). 
Since $\eE[\psi_\tau(Y_i-g(\X_i; \eb^1))]=0$ for any $i \geq m+k^0_m+1$, we have:
\begin{equation}
\label{eq13bis}
  \eE[\sum^{m+\widetilde k}_{i=m+k^0_m+1} \eR_i(\UU) ] =\sum^{m+\widetilde k}_{i=m+k^0_m+1}   \eg(\X_i;\ebo+m^{-1/2} \UU) \big[ F(0)   -F\big(g(\X_i;\eb^1)- g(\X_i; \ebo+m^{-1/2} \UU)\big) \big] .
  \end{equation}
  For continuing   the study of $\sum^{m+\widetilde k}_{i=m+k^0_m+1} \eR_i(\UU) $, we consider the two possible cases for difference $\eb^1 -\ebo$. \\
\hh \textit{(a)}  If  $\| \eb^1-\ebo \|>c_1>0$, then by a Taylor expansion of the distribution function $F$ in relation (\ref{eq13bis}), we have:
  \[
      \eE[\sum^{m+\widetilde k}_{i=m+k^0_m+1} \eR_i(\UU) ]  = \sum^{m+\widetilde k}_{i=m+k^0_m+1}   \eg(\X_i;\ebo+m^{-1/2} \UU)   \big[ g(\X_i; \ebo+m^{-1/2} \UU))- g(\X_i; \eb^1)\big] f(\widetilde{\widetilde{b}_i}) ,    \]
  with  the    $\widetilde{\widetilde{b}_i}$ between $0$ and $g(\X_i;\eb^1)- g(\X_i; \ebo+m^{-1/2} \UU)$.\\
 Since there exists a constant $c_3>0$ such that $f(\widetilde{\widetilde{b}_i}) >c_3$ and since, by assumption (A8), 
 \[
 \frac{1}{\widetilde k - k^0_m} \left\| \sum^{m+\widetilde k}_{i=m+k^0_m+1} \eg(\X_i;\ebo) \big[g(\X_i; \eb^1) - g(\X_i; \eb^0) \big] \right\|_{\infty} >c_2>0,
 \]
then, we have:
    \[
  \eE[\|\sum^{m+\widetilde k}_{i=m+k^0_m+1} \eR_i(\UU) \| ]= O(\widetilde k - k^0_m) .
  \]
  By assumption (A3) and the strong law of large numbers, together with the last relation, we have that
   \begin{equation}
    \label{ED0}
\big\| \sum^{m+\widetilde k}_{i=m+k^0_m+1} \eR_i(\UU) \big\| =O_{\PP} (\eE[\|\sum^{m+\widetilde k}_{i=m+k^0_m+1} \eR_i(\UU) \| ])= O_{\PP} (\widetilde k - k^0_m) =O_{\PP}(m^{s-1/2}\| m^{1/2}(\eb^1-\ebo)\|).
  \end{equation}
\hh \textit{(b)}  Consider now the case  $\eb^1-\ebo \rightarrow \textbf{0}$, such that $m^{1/2} \| \eb^1- \ebo \| \rightarrow\infty$,  as $m \rightarrow \infty$.\\
  Since $\eg(\xx,\eb)$ is bounded for any $\xx \in \Upsilon$ and for any $\eb$ such that $\| \eb -\ebo \| \leq b_m$, we have:
  \[
 \max_{m+k^0_m+1 \leq i \leq m+\widetilde k} \bigg( \max_{\substack{\| \UU\| \leq C_1\\\| \eb - \ebo\| \leq \| \eb^1- \ebo\|}} \| \eR_i(\UU) \|\bigg)= O_{\PP}(1).
 \]
 In this case, relation (\ref{AB}) becomes:
 $R_{ij}(\UU)=A_{ij}(\UU)+B_{ij}(\UU)$, with
 \[
 A_{ij}(\UU)= \frac{\partial g((\X_i;\ebo+m^{-1/2} \UU))}{\partial \beta_j}\big[\e1_{\varepsilon_i \leq 0} -\e1_{\varepsilon_i \leq g(\X_i; \ebo+m^{-1/2} \UU)-g(\X_i; \eb^1)} \big],
 \]
 \[  B_{ij}(\UU)=\big[ \frac{\partial g(\X_i;\ebo+m^{-1/2} \UU)}{\partial \beta_j}  -\frac{\partial g(\X_i;\ebo)}{\partial \beta_j}\big][\tau - \e1_{\varepsilon_i \leq 0}].
 \]
 Similarly as for relation (\ref{eq20}), we have uniformly in $\UU$:
 $
 \sum^{m+\widetilde k}_{i=m+k^0_m+1} \eE[B_{ij}^2(\UU)] =O((\widetilde k- k^0_m)m^{-1})$.\\
 For random variable $A_{ij}(\UU)$, we have similarly as in the proof of Lemma \ref{Lemma 1}, taking into account that $m^{1/2}\|\eb^1-\ebo\| \rightarrow \infty$, and that $\eg$ is bounded in the neighbourhood of $\ebo$, that: $ \sum^{m+\widetilde k}_{i=m+k^0_m+1} \eE[A_{ij}^2(\UU)]  = O((\widetilde k- k^0_m)\| \eb^1 - \ebo\|)$.\\
 Then 
 \[
 \sum^{m+\widetilde k}_{i=m+k^0_m+1} Var[R_{ij}(\UU)] \leq  O((\widetilde k- k^0_m)\| \eb^1 - \ebo\|).
 \]
 We take $\delta_m=(\widetilde k- k^0_m)\| \eb^1 - \ebo\|$ and applying Lemma \ref{Lemma 4.4.1.}, we get:
 \begin{equation}
 \label{DED}
 \sum^{m+\widetilde k}_{i=m+k^0_m+1} \eR_i(\UU)=\eE[\sum^{m+\widetilde k}_{i=m+k^0_m+1} \eR_i(\UU)]+O_{\PP}(m^{s/2}\sqrt{s \log m} \sqrt{\| \eb^1 - \ebo\|}).
 \end{equation}
 We study now $  \eE[\sum^{m+\widetilde k}_{i=m+k^0_m+1} \eR_i(\UU)] $. Taking into account relation (\ref{eq10}), by a Taylor expansion, we get
   \[
      \eE[\sum^{m+\widetilde k}_{i=m+k^0_m+1} \eR_i(\UU)] = \sum^{m+\widetilde k}_{i=m+k^0_m+1}   \eg(\X_i;\ebo+m^{-1/2} \UU) \bigg[(-m^{-1/2} \UU+\eb^1-\ebo)^t\eg(\X_i;\ebo+m^{-1/2} \UU )
      \]
      \[
   \qquad  +\frac{(\eb^1-\ebo)^t}{2}\egg(\X_i; \widetilde{\eb})(\eb^1-\ebo) \bigg] f(\widetilde{\widetilde{b}_i}) .
     \]
Since $m^{1/2}\|\eb^1-\ebo\| \rightarrow \infty$ as $m \rightarrow \infty$,    taking into account assumptions (A2), (A4) and (A5), we have, uniformly in $\UU$ bounded:
  \[
  \begin{array}{cl} 
 \displaystyle{  \eE[\sum^{m+\widetilde k}_{i=m+k^0_m+1} \eR_i(\UU)] }& = (\eb^1-\ebo)^t \displaystyle{\sum^{m+\widetilde k}_{i=m+k^0_m+1}  } \eg(\X_i;\ebo+m^{-1/2} \UU) \bigg\{ \eg(\X_i;\ebo+m^{-1/2} \UU)  f(\widetilde{\widetilde{b}_i})  \bigg\} (1+o(1))\\
    &= O(m^{s}(\eb^1-\ebo)).
       \end{array}
 \]
 Hence,  replacing this in (\ref{DED}): 
   \begin{equation}
    \label{ED1}
  \sum^{m+\widetilde k}_{i=m+k^0_m+1} \eR_i(\UU)=  O(m^{s}(\eb^1-\ebo))+ O_{\PP}(m^{s/2}\sqrt{s \log m} \sqrt{\| \eb^1 - \ebo\|})=O_{\PP}(m^{s}(\eb^1-\ebo)).
     \end{equation}
  Thus, in the two cases of  $\eb^1-\ebo$, taking into account (\ref{ED0}) and (\ref{ED1}),    we have
  \begin{equation}
  \label{eq11}
   \|\sum^{m+\widetilde k}_{i=m+k^0_m+1} \eR_i(\UU)\| = O_{\PP}(m^{s-1/2 }\| m^{1/2}(  \eb^1-\ebo)\|).
  \end{equation}
  On the other hand, for (\ref{eq10}), we have by functional central limit theorem:
  \begin{equation}
  \label{eq12}
  \sum^{m+\widetilde k}_{i=m+k^0_m+1}  \eg(\X_i;\eb^0)\psi_\tau(\varepsilon_i) =O_{\PP}(m^{s/2}).
  \end{equation}
  Taking into account  relations  (\ref{eq10}),   (\ref{eq11}), (\ref{eq12}), and $s>1$, we obtain:
  \[
\big\|  \J_m^{-1/2}(\widehat \eb_m)\sum^{m+\widetilde k}_{i=m+k^0_m+1} \eg(\X_i;\ebo+m^{-1/2} \UU)  \psi_\tau\big(Y_i-g(\X_i; \ebo+m^{-1/2} \UU)\big)\big\|_{\infty}/z(m,\widetilde k,\gamma)
  \]
 \begin{equation}
  \label{eq14}
  = \frac{O_{\PP}(m^{s/2})+O_{\PP}(m^{s-1/2}\|m^{1/2} ( \eb^1-\ebo)\|_{\infty})}{\sqrt{m} (1+\widetilde k/m)(\widetilde k/(m+\widetilde k))^\gamma}  .
  \end{equation}
Since $\widetilde k =k^0_m+m^s$, we have  $ \lim_{m \rightarrow \infty}\left(\frac{\widetilde k}{m+\widetilde k} \right)^\gamma = 1$  and $\sqrt{m} \left( 1+\frac{\widetilde k}{m}\right) =m^{s-1/2}$. Then, relation (\ref{eq14}) converges  to infinity as  $m \rightarrow \infty$.  
Thus, 
 \[
  \frac{\left\| \J_m^{-1/2}(\widehat \eb_m)\sum^{m+\widetilde k}_{i=m+k^0_m+1} \eg(\X_i;\widehat \eb_m)  \psi_\tau\big(Y_i-g(\X_i; \widehat \eb_m)\big) \right\|_{\infty}}{z(m,\widetilde k,\gamma)} \overset{{\PP}} {\underset{m \rightarrow \infty}{\longrightarrow}}  \infty .
  \]
 Theorem follows for the open-end procedure taking  into account the last relation together with relations  (\ref{eq9}), (\ref{TJ}).\\ 
  
{\textit{(ii) Closed-end procedure case.}}  The proof is similar to that of \textit{(i)}, taking   $\widetilde k =k^0_m+m$. In this case, by analogous calculations,   relation (\ref{eq11}) becomes $\|\sum^{m+\widetilde k}_{i=m+k^0_m+1} \eR_i(\UU)\| \geq O_{\PP}(m^{1/2}\|m^{1/2}(\eb^1-\ebo)\|)$ and relation (\ref{eq12}) becomes $ \sum^{m+\widetilde k}_{i=m+k^0_m+1}  \eg(\X_i;\ebo)\psi_\tau(\varepsilon_i) =O_{\PP}(m^{1/2})$. Therefore,  relation (\ref{eq14}) converges  to infinity as  $m \rightarrow \infty$.  
 \hspace*{\fill}$\blacksquare$   \\


\end{document}